\documentclass [a4paper] {article}
\usepackage {natbib}

\usepackage {amsmath}
\usepackage {amsfonts}
\usepackage {amsthm}
\usepackage{hyperref}
\makeatletter
\newcommand \th@slplain {\th@plain \slshape}
\makeatother
\theoremstyle {slplain}
\newtheorem {theorem}{Theorem}
\newtheorem {corollary}[theorem]{Corollary}
\theoremstyle {remark}
\newtheorem {example}[theorem]{Example}
\newtheorem {remark}[theorem]{Remark}

\makeatletter
\newcommand* \eqarr@ywithdef [1]{%
	\null\,%
	\vcenter {%
		\normalbaselines
		\openup\jot
		\m@th
		\ialign {\strut #1\crcr }%
	}\,%
}
\newcommand* \eqalign [1]{\eqarr@ywithdef {%
	\hfil$\displaystyle{##}$&$\displaystyle{{}##}$\hfil \crcr
	#1%
}}
\newcommand* \listlines [1]{\eqarr@ywithdef {%
	\hfil$\displaystyle{##}$\hfil \crcr
	#1%
}}
\newcommand* \eqarraywithdef [1]{%
	\eqarr@ywithdef {\ignorespaces #1}%
}
\makeatother

\newcommand* \ba {\mathop{\rm ba}\nolimits }
\newcommand* \core {\mathop{\rm core}\nolimits }
\newcommand* \field {\mathop{\rm field}\nolimits }
\newcommand* \bacore {\mathop {\rm ba\mathchar`-core}\nolimits}

\begin {document}

\title {\kern-10pt
	On Balanced Games with Infinitely Many Players:\kern-10pt\\
	Revisiting Schmeidler's Result}

\author {%
	\llap {\vtop {\hsize=200pt \noindent
		David Bartl\endgraf
		\medskip\smallskip
		\vtop {\small \noindent
			Department of Informatics and Mathematics,\break
			School of Business Administration in Karvin\'a,\break
			University in Opava\break
			{\tt bartl@opf.slu.cz}}}}%
	\rlap {\vtop {\hsize=200pt \noindent
		Mikl\'os Pint\'er\endgraf
		\medskip\smallskip
		\vtop {\small \noindent
			Corvinus Center for Operational Research,\break
			Institute of Advanced Studies,\break
			Corvinus University of Budapest\break
			{\tt pmiklos@protonmail.com}}}}}




\begingroup
\let \newpageori = \newpage
\def \newpage {\newpageori \null \vskip-60pt \let\newpage=\empty}
\maketitle
\let \newpage = \newpageori
\endgroup


\begin{abstract}\noindent
We consider transferable utility cooperative games with infinitely many
players and the core understood in the space of bounded additive set functions.  We
show that, if a game is bounded below, then its core is non-empty if
and only if the game is balanced.  

This finding is a generalization of Schmeidler's (1967) original result
``On Balanced Games with Infinitely Many Players'', where the game is
assumed to be non-negative.  We furthermore demonstrate that, if a game is
not bounded below, then its core might be empty even though the game is
balanced; that is, our result is tight. 

We also generalize Schmeidler's (1967) result to the case of restricted cooperation too.
\end {abstract}

\begingroup

\begin {abstract}\noindent
	TU games with infinitely many players, core, balancedness, TU games with restricted cooperation, signed TU games, bounded additive set functions
\end {abstract}
\endgroup

\quad \noindent {\it 2020 Mathematics Subject Classification:\/}
	91A12, 91A07

\quad \noindent {\it JEL Classification:\/}
	C71

\section {Introduction}

The core \citep{Shapley1955,Gillies1959} is one of the most important solution concepts of cooperative game theory. It is important not only from the theory viewpoint, but for its simple and easy
to understand nature, it also helps to solve various problems arising in practice. 

In the transferable utility setting (henceforth TU games) the Bondareva-Shapley Theorem \citep{Bondareva1963,Shapley1967,Faigle1989} provides a necessary and sufficient condition for the non-emptiness of the core of a finite TU game; it states that the core of a finite TU game with our without restricted cooperation is not empty if and only if the TU game is balanced. The textbook proof of the Bondareva-Shapley Theorem goes by the strong duality theorem of linear programs, see e.g.\ \citet{PelegSudholter2007}. 

\citet{Schmeidler1967}, \citet{Kannai1969,Kannai1992}, and \citet{Pinter2011d}, among others, considered TU games with infinitely many players. All these papers studied the case when the core consists of bounded additive set functions. \citet{Schmeidler1967} and \citet{Kannai1969} showed respectively that the core of a non-negative TU game with infinitely and countably infinitely many players is not empty if and only if the TU game is balanced.

In this paper we consider infinite signed TU games (sign unrestricted TU games with infinite many players) with and without restricted cooperation. Particularly, we follow \citet{Schmeidler1967} and assume that the allocations are bounded additive set functions. 

Applications of infinite signed TU games go back in times at least as early as \citet{ShapleyShubik1969b} (economic systems with externalities), which generalize market games \citep{ShapleyShubik1969}. Further applications are (semi-) infinite transportation games \citep{Sanchez-SorianoLlorcaTijsTimmer2002,TimmerLlorca2002}, infinite sequencing games \citep{FragnelliLlorcaSanchezSorianoTijsBranzei2010}, and somehow less directly the line of literature represented by  e.g. \citet{MontrucchioSemeraro2008} among others. 

While we can analyze the non-emptiness of the core in the finite setting by using the aforementioned Bondareva-Shapley Theorem \citep{Bondareva1963,Shapley1967,Faigle1989}, we have been missing an appropriate tool for such TU games with infinitely many players.

Our contribution is an extension of Schmeidler's \citeyearpar{Schmeidler1967} result saying a non-negative infinite TU game without restricted cooperation has a non-empty core if and only if it is balanced, to the general case saying a bounded below infinite TU game with or without restricted cooperation has a non-empty core if and only if it is balanced (Theorems \ref{thm-Bondareva-Shapley-megszamjatek-Schmeidler} and \ref{megszamjatek-Schmeidler-proposition}).

It is worth mentioning that neither Schmeidler's \citeyearpar{Schmeidler1967} nor Kannai's \citeyearpar{Kannai1969,Kannai1992} approach (proof) can be applied to achieve our generalization (Theorems \ref{thm-Bondareva-Shapley-megszamjatek-Schmeidler} and \ref{megszamjatek-Schmeidler-proposition}). Our approach is different from the previous ones.


The set-up of this paper is as follows.  In Sections
\ref{sec1}~and~\ref{sec2}, we introduce basic notions of TU games with
infinitely many players, including the core and balancedness, and we
present our main result (Theorem~\ref
{thm-Bondareva-Shapley-megszamjatek-Schmeidler}).  In Sections
\ref{sec3}~and~\ref{sec4}, we recall some useful concepts pertaining
functional spaces, topology and compactness, and we prove our main
result.  We additionally give examples to show the tightness of our main
result and we also mention an interesting ``limiting'' property of the
core.  Finally, in Section~\ref{sec6}, we discuss the case of restricted
cooperation and give our second main result (Theorem~\ref
{megszamjatek-Schmeidler-proposition}).

\section{Preliminaries of infinite TU games}\label{sec1}

We consider transferable utility cooperative games with a finite or
infinite set~$N$ of players.  A {\it coalition\/} is a subset $S
\subseteq N$, so the power set $\mathcal P(N) = \{\, S : S \subseteq N
\,\}$ is the collection of all coalitions that can be considered.
Let $\mathcal A \subseteq \mathcal P(N)$ be the collection of all
\textit{feasible coalitions}, which are those that can potentially emerge.  In the case of no restricted cooperation, we assume that $\mathcal A$~is
a field of sets over~$N$; that is, the collection $\mathcal A$~is such
that $\emptyset \in \mathcal A$ and, if $S, T \in \mathcal A$, then $N
\setminus S \in \mathcal A$ and also $S \cup T \in \mathcal A$. In the case of restricted cooperation, we assume only that $\emptyset,N \in \mathcal{A}$.

Then a {\it transferable utility cooperative game\/}
(henceforth {\it game\/} for short) is represented by its {\it
coalition function}, which is a mapping $v\colon \mathcal A \to \mathbb
R$ such that $v(\emptyset) = 0$.  For any coalition $S \in \mathcal A$,
the value $v(S)$ is understood as the payoff that the coalition
$S$~receives if it is formed.

Assume that the players form the grand coalition $N \in \mathcal A$.
Then $v(N) \in \mathbb R$ is the value of the grand coalition~$N$, and
the issue is to allocate this value among the players.
Following \citet {Schmeidler1967}, we define the allocations
as bounded additive set functions $\mu\colon \mathcal A
\to \mathbb R$; that is, a function such that $\bigl|\mu(S)\bigr| \leq
C$ for all $S \in \mathcal A$ for some constant $C \in \mathbb
R$ and $\mu(S \cup T) = \mu(S) + \mu(T)$ for any disjoint $S,T \in
\mathcal A$.  Let $\ba(\mathcal A) = \{\, \mu\colon \mathcal A \to
\mathbb R : \mu$~is a bounded additive set function$\,\}$ denote the
space of all bounded additive set functions on~$\mathcal A$.  Then the
{\it core\/} of the game by $v$ is the set
$$
\eqarraywithdef {\hfil $# \mu($& $\hfil #\hfil$& $)\hfil #\hfil
	v($& $\hfil #\hfil$& $)#\hfil $\cr
\bacore(v) = \bigl\{\, \mu \in \ba(\mathcal A) : & N & = & N & \,, \cr
 & S & \geq & S & \qquad \hbox {for all} \quad S \in \mathcal A \setminus \{N\}
		\,\bigr\} \,\hbox{.} }
$$
In words, the core consists of all the allocations of the value
$v(N)$ among the players (efficiency) such that any
coalition $S \in \mathcal A \setminus \{N\}$ that could potentially
emerge gets by the proposed allocations at least as much as the value $v(S)$ (coalitional rationality), see \citet {Shapley1955}, \citet
{Gillies1959}, \citet {Kannai1992}, and \citet {Zhao2018}. 

It is worth noticing that calling any game where the class of feasible coalitions is a field even if the field is not the power set of the player set is not misleading because any additive set function defined on a subfield can be extended to the power set. Therefore, an allocation from the core of a game without restricted cooperation gives rise (typically in a non-unique way) to an allocation defined on the power set of the player set. 

The case when the class of feasible coalitions is not a field, however, leads to the very same features of the core as restricted cooperation leads in the finite setting (see \citet{Faigle1989}), explaining why we call this case restricted cooperation.

The key question is whether the core is non-empty. An answer is provided by the Bondareva-Shapley Theorem.

\section {The Bondareva-Shapley Theorem}\label {sec2}

Consider a game having finitely many players without restricted cooperation.  In this case,
we have $N = \{1,2,\ldots,n\}$ for some natural number~$n$ and $\mathcal
A = \mathcal P(N)$.  Moreover, in this setting, the allocations of the
value $v(N)$ among the players are given by \textit{payoff vector}s,
any of them is an $n$-tuple $a = (a_i)_{i=1}^n \in \mathbb R^N$ of real
numbers; number $a_i$~means the payoff allocated to player~$i$ for
$i = 1$,~$2$,~\dots,~$n$.  Then the {\it core\/} of this game is
defined to be the set
$$
\eqarraywithdef {\hfil $# {}$& $\mathrel {} \sum_{i\in #} \hfil a_i$&
	${}\hfil #\hfil v($& $\hfil #\hfil$& $)#\hfil $\cr
\core(v) = \bigl\{\, a \in \mathbb R^N : & N & = & N & \,, \cr
 & S & \geq & S & \qquad \hbox {for all} \quad
			S \in \mathcal P(N) \setminus \{N\}
		\,\bigr\} \,\hbox{.} }
$$
The intuitive meaning of the $\core(v)$ is the same as that of
the $\bacore(v)$, see above.  Clearly, given a payoff vector $a \in \mathbb
R^N$, we can define the corresponding additive set function $\mu\colon
\mathcal P(N) \to \mathbb R$ by $\mu(S) = \sum_{i\in S} a_i$ for any $S
\in \mathcal P(N)$.  Conversely, given an additive set function
$\mu\colon \mathcal P(N) \to \mathbb R$, we can define the corresponding
payoff vector $a \in \mathbb R^N$ by $a_i = \mu\bigl(\{i\}\bigr)$ for $i
= 1$,~$2$,~\dots,~$n$.  Here any additive set function $\mu\colon
\mathcal P(N) \to \mathbb R$ is bounded as the number of the players is
finite.  We thus have a one-to-one correspondence between the
$\bacore(v)$ and the $\core(v)$.  Hence, the notion of bounded additive
function $\mu \in \ba(\mathcal A)$ naturally extends the concept of the
payoff vector $a \in \mathbb R^N$ when the set~$N$ of the players is
infinite.

Regarding the question whether the $\core(v)$ is non-empty, for any
coalition $S \subseteq N$, define its {\it characteristic vector\/} to
be the row vector $\chi_S = \bigl( \left.\chi_S(1)\right. \break\,
\left.\chi_S(2)\right. \, \left.\ldots\right. \allowbreak\,
\left.\chi_S(n)\right.  \bigr)$ with $\chi_S(i) = 1$ if $i \in S$, and
with $\chi_S(i) = 0$ if $i \notin S$, for $i = 1$,~$2$,~\dots,~$n$.  We
say that a collection $\mathcal S = \{S_1, S_2, \ldots, \allowbreak
S_r\} \subseteq \mathcal P(N)$ of coalitions is {\it balanced\/} if
there exist non-negative real numbers
$\lambda_1$,~$\lambda_2$,~\dots,~$\lambda_r$, called balancing
weights, such that
\begin {equation}
\label {eq-balanced-set-of-coalitions}
\sum_{p=1}^r \lambda_p \chi_{S_p} = \chi_N \,\hbox{.}
\end {equation}
Moreover, we say that a game $v$ is {\it balanced\/} if
\begin {equation}
\label {eq-balanced-game}
\sum_{p=1}^r \lambda_p v(S_p) \leq v(N)
\end {equation}
for every balanced collection $\{S_1, S_2, \ldots, \allowbreak S_r\}
\subseteq \mathcal P(N)$ of coalitions.  The following result due to
\citet {Bondareva1963} and \citet {Shapley1967}, later extended by \citet{Faigle1989} to the restricted cooperation case, has become classical:

\begin {theorem} [Bondareva-Shapley Theorem]
\label {thm-Bondareva-Shapley}
Consider a game with finitely many players, with or without restricted cooperation, represented by a
coalition function\/ $v\colon \mathcal P(N) \to \mathbb R$.  Then the\/
$\core(v)$ is non-empty if and only if the game is balanced.
\end {theorem}

Consider now a general game without restricted cooperation; that is, the set~$N$ of the players can be
finite or infinite and  the class of feasible coalitions $\mathcal A \subseteq \mathcal P(N)$ is a field of sets over~$N$.  Concerning the question
whether $\bacore(v)$ is non-empty we follow \citet
{Schmeidler1967}, who proceeds analogously as in the classical case;
that is:

For any subset $S \subseteq N$, define its {\it characteristic function\/}
$\chi_S\colon N \to \{0,1\}$ by letting $\chi_S(i) = 1$ if $i \in S$,
and $\chi_S(i) = 0$ if $i \notin S$, for every $i \in N$.  We say that a
collection $\mathcal S = \{S_1, S_2, \ldots, \allowbreak S_r\} \subseteq
\mathcal A$ of coalitions is {\it balanced\/} if there exist
non-negative real numbers $\lambda_1$,~$\lambda_2$,~\dots,~$\lambda_r$,
called balancing weights, such that
\begin {equation}
\label {eq-Schmeidler-balanced-set-of-coalitions}
\sum_{p=1}^r \lambda_p \chi_{S_p} = \chi_N \,\hbox{.}
\end {equation}
Furthermore, we say that a game $v$ is {\it balanced\/} if
\begin {equation}
\label {eq-Schmeidler-balanced-game}
\sum_{p=1}^r \lambda_p v(S_p) \leq v(N)
\end {equation}
for every balanced collection $\{S_1, S_2, \ldots, \allowbreak S_r\}
\subseteq \mathcal A$ of coalitions.

\begin {remark}
\citet {Schmeidler1967} actually defines balancedness in a slightly
different way:  ``A game is balanced if $\sup \sum_i a_i v(A_i)
\leq v(S)$ when the sup is taken over all finite sequences of
$a_i$~and~$A_i$, where the $a_i$~are non-negative numbers,
the $A_i$~are in~$\Sigma$, and $\sum_{i} a_i \chi_{A_i}
\leq \chi_S$.'' Considering non-negative games, \citeauthor {Schmeidler1967} explains that his definition is different from the ``definition with equality" only in its form:  ``It is easy to verify that this sup does not change
even if it is constrained by $\sum_i a_i \chi_{A_i} =
\chi_S$ (instead of the inequality); also, for balanced games, the sup
equals $v(S)$.''~--- See \citet [p.~1] {Schmeidler1967}. In the case of non-negative games \citeauthor{Schmeidler1967}'s definition of balancedness is equivalent with the ``definition with equality"; however, in the general, signed case, those are different.
\end {remark}

Then \citet{Schmeidler1967} proves the following result, see \citet
{Kannai1969} for another proof:

\begin {theorem} [Bondareva-Shapley Theorem, \citet {Schmeidler1967}]
\label {thm-Bondareva-Shapley-Schmeidler}
Given a finite or infinite set\/~$N$ of the players and a field of
sets\/ $\mathcal A \subseteq \mathcal P(N)$ over\/~$N$, consider a
game represented by a coalition function\/ $v\colon \mathcal A \to
\mathbb R$.  If the game is non-negative; that is,
$$
\forall S \in \mathcal A\colon\,\,\,
v(S) \geq 0 \,,
$$
then the\/ $\bacore(v)$ is non-empty if and only if the game is balanced.
\end {theorem}

It is easy to see that Theorem~\ref {thm-Bondareva-Shapley-Schmeidler}
is a generalization of Theorem~\ref {thm-Bondareva-Shapley} if the game
is non-negative.  Our goal, nonetheless, is to establish the following
result:

\begin {theorem} [Bondareva-Shapley Theorem, a generalization]
\label {thm-Bondareva-Shapley-megszamjatek-Schmeidler}
Given a finite or infinite set\/~$N$ of the players and a field of
sets\/ $\mathcal A \subseteq \mathcal P(N)$ over\/~$N$, consider a
game represented by a coalition function\/ $v\colon \mathcal A \to
\mathbb R$.  If the game is bounded below; that is,
$$
\exists L \in \mathbb R \,\,\,\,
\forall S \in \mathcal A\colon\,\,\,
v(S) \geq L \,,
$$
then the\/ $\bacore(v)$ is non-empty if and only if the game is balanced.
\end {theorem}

Notice that Theorem~\ref {thm-Bondareva-Shapley-megszamjatek-Schmeidler}
directly generalizes both Theorems \ref {thm-Bondareva-Shapley}~and~\ref
{thm-Bondareva-Shapley-Schmeidler}
because a game with finitely many players is always bounded below.

Before we present our proof of Theorem~\ref
{thm-Bondareva-Shapley-megszamjatek-Schmeidler}, we find it appropriate
to introduce and recall several notions and concepts.

\section {Several notions and concepts}\label {sec3}

Let $N$~be a set and let $\mathcal A \subseteq \mathcal P(N)$ be a field
of sets over~$N$.  Then the pair $(N, \mathcal A)$ is called {\it
chargeable space}.  Recall that, for any $S \subseteq N$, the symbol
$\chi_S$~denotes the characteristic function $\chi_S\colon N \to
\{0,1\}$ of the set~$S$.  Given a function $f\colon N \to \mathbb R$, we
say it is a {\it simple function\/} if $f = \lambda_1
\chi_{S_1} + \lambda_2 \chi_{S_2} + \cdots + \lambda_r \chi_{S_r}$ for
some natural number~$r$, for some real numbers
$\lambda_1$,~$\lambda_2$,~\dots,~$\lambda_r$, and for some sets $S_1,
S_2, \ldots, \allowbreak S_r \in \mathcal A$. Let $\Lambda(\mathcal A) = \{\, f\colon N \to \mathbb R : f$~is a simple function$\,\}$ denote the vector (i.e.\ linear) space of all
simple functions defined over $(N,\mathcal{A})$, where the sum of two
functions and the multiplication of a function by a constant are both
defined in the usual way, i.e.\ pointwise.  For a simple function $f \in \Lambda(\mathcal A)$, define its norm to be
$$
\mathopen\|f\mathclose\| = \sup_{i\in N} \bigl|f(i)\bigr| \,\hbox{,}
$$
so $\Lambda(\mathcal A)$ is a normed linear space.  

Likewise, notice that the space $\ba(\mathcal A)$ of all bounded
additive set functions on~$\mathcal A$ is also a vector space; for a $\mu
\in \ba(\mathcal A)$, define its norm to be
\begin {equation}
\label {eq-mu-norm}
\mathopen\|\mu\mathclose\| = \sup_{\substack{
	r \in \mathbb N \\
	S_1, S_2, \ldots, S_r \in \mathcal A \\
	\hidewidth S_1 \cup S_2 \cup \cdots \cup S_r = N \hidewidth \\
	S_i \cap S_j = \emptyset, \; i \neq j
}} \bigl|\mu(S_1)\bigr| + \bigl|\mu(S_2)\bigr| + \cdots + \bigl|\mu(S_r)\bigr|
\,\hbox{.}
\end{equation}

It is well-known that the topological dual $(\Lambda(\mathcal
A))^*$ of the vector space $\Lambda(\mathcal A)$; that is, the space of
all continuous linear functionals on $\Lambda(\mathcal A)$, is
isometrically isomorphic to the space $\ba(\mathcal A)$ (see, e.g., \citet{DunfordSchwartz1958}, Theorem IV.5.1, p.~258).  Indeed, a
continuous linear functional $\mu' \in (\Lambda(\mathcal A))^*$ induces
a bounded additive set function $\mu \in \ba(\mathcal A)$ by letting
$\mu(S) = \mu'(\chi_S)$ for $S \in \mathcal A$, and, conversely, a
bounded additive set function $\mu \in \ba(\mathcal A)$ induces a
continuous linear functional $\mu' \in (\Lambda(\mathcal A))^*$ by
letting
\begin {equation}
\label {eq-funct}
\mu'(f) = \lambda_1 \mu(S_1) + \lambda_2 \mu(S_2) + \cdots + \lambda_r \mu(S_r)
\end{equation}
for any simple function $f = \lambda_1 \chi_{S_1} + \lambda_2 \chi_{S_2}
+ \cdots + \lambda_r \chi_{S_r} \!\in \Lambda(\mathcal A)$.  This is the
reason why, for simplicity, we shall identify the space
$(\Lambda(\mathcal A))^*$ with $\ba(\mathcal A)$.

Consider now a game represented by a coalition function $v\colon
\mathcal A \to \mathbb R$, and let the game be bounded below; that
is, there exists a constant $L \in \mathbb R$ such that $v(S) \geq L$
for all $S \in \mathcal A$.  Assume that a $\mu \in \bacore(v)$.  Let
$S_1, S_2, \ldots, \allowbreak S_r \in \mathcal A$ be pairwise disjoint
and such that $N = S_1 \cup S_2 \cup \cdots \cup S_r$.  Then
$$
\eqarraywithdef {\hfil $\displaystyle #$& $\displaystyle {}#$\hfil \cr
\sum_{p=1}^r \bigl|\mu(S_p)\bigr| &=
\sum_{\substack{ p=1 \\ \hidewidth \mu(S_p) \geq 0 \hidewidth}}^r \mu(S_p) -
\sum_{\substack{ p=1 \\ \hidewidth \mu(S_p)   <  0 \hidewidth}}^r \mu(S_p) \cr &=
\mu\biggl(\bigcup_{\substack{ p=1 \\ \hidewidth \mu(S_p)\geq0 \hidewidth}}^r S_p\biggr) -
\mu\biggl(\bigcup_{\substack{ p=1 \\ \hidewidth \mu(S_p)  < 0 \hidewidth}}^r S_p\biggr) \cr &=
\mu\biggl(N \setminus \bigcap\limits_{\substack{ p=1 \\ \hidewidth \mu(S_p) \geq 0 \hidewidth}}^r (N \setminus S_p)\biggr) -
\mu\biggl(\bigcup_{\substack{ p=1 \\ \hidewidth \mu(S_p)  < 0 \hidewidth}}^r S_p\biggr) \cr &=
\mu(N) -
\mu\biggl(\bigcap_{\substack{ p=1 \\ \hidewidth \mu(S_p)\geq0 \hidewidth}}^r (N \setminus S_p)\biggr) -
\mu\biggl(\bigcup_{\substack{ p=1 \\ \hidewidth \mu(S_p)  < 0 \hidewidth}}^r S_p\biggr) \cr &\leq
\mu(N) - 2L = v(N) - 2L \,\hbox{.} }
$$
By taking the definition~\eqref {eq-mu-norm} of the norm into account,
it follows the $\bacore(v)$ is contained in the closed ball $B_R =
\bigl\{\, \mu \in \ba(\mathcal A) : \mathopen\|\mu\| \leq v(N) - 2L
\,\bigr\}$ of radius $R = v(N) - 2L$.  (Notice that $v(N) - 2L \geq 0$,
for if we had $v(N) < 2L$, then the $\bacore(v)$ would obviously be
empty, contradicting the assumption that $\mu \in \bacore(v)$.)

We endow the space $\ba(\mathcal A)$ with the weak* topology with
respect to $\Lambda(\mathcal A)$.  The topology will be introduced if we
describe all the neighborhoods of a point.  A set $U \subseteq
\ba(\mathcal A)$ is a {\it weak* neighborhood\/} of a $\mu_0 \in
\ba(\mathcal A)$ if there exist a natural number~$r$ and functions
$f_1, f_2, \ldots, \allowbreak f_r \in \Lambda(\mathcal A)$ such that
$\bigcap_{p=1}^r \bigl\{\, \mu \in \ba(\mathcal A) : \bigl|\mu'(f_p) -
\mu'_0(f_p)\bigr| < 1 \,\bigr\} \subseteq U$, where $\mu'$~and~$\mu'_0$
is the continuous linear functional induced by $\mu$~and~$\mu_0$,
respectively, see~\eqref{eq-funct}.  By Alaoglu's Theorem (see, e.g.,
\citet {AliprantisBorder2006}, Theorem~6.21, p.~235), the closed ball
$B_R$~is compact in the weak* topology.  That is, if $G_i \subseteq
\ba(\mathcal A)$ are weakly* open sets for $i \in I$, where $I$~is an
index set, such that $\bigcup_{i\in I} G_i \supseteq B_R$, then
$\bigcup_{j=1}^n G_{i_j} \supseteq B_R$ for some natural number~$n$ and
for some $i_1, i_2, \ldots, \allowbreak i_n \in I$.

Let $F_i \subseteq B_R$ be weakly* closed sets for $i \in I$, where
$I$~is an index set.  We say the collection $\{F_i\}_{i\in I}$ is a {\it
centered system\/} of sets if $\bigcap_{j=1}^n F_{i_j} \neq \emptyset$
for any natural number~$n$ and for any $i_1, i_2, \ldots, \allowbreak
i_n \in I$.  By considering the complements ($G_i = \ba(\mathcal A)
\setminus F_i$), it follows $\bigcap_{i\in I} F_i \neq \emptyset$.

In our proof of
Theorem~\ref{thm-Bondareva-Shapley-megszamjatek-Schmeidler} we consider
the weakly* closed sets
$$
F_S = \bigl\{\, \mu \in \ba(\mathcal A) :
	\mu(N) = v(N) \text{ and }
	\mu(S) \geq v(S) \text{ and }
	\mathopen\|\mu\mathclose\| \leq R \,\bigr\}
$$
for $S \in \mathcal A$.  The main idea is to show that, if the game
$v$~is balanced, then the system $\{F_S\}_{S\in\mathcal A}$
is centered.  Noticing that $\bacore(v) = \bigcap_{S\in\mathcal A} F_S
\neq \emptyset$, the proof will be done.

We are now ready to present our proof of Theorem~\ref{thm-Bondareva-Shapley-megszamjatek-Schmeidler}.

\section {Proof of Theorem~\ref{thm-Bondareva-Shapley-megszamjatek-Schmeidler}}\label {sec4}


Below we give our proof of Theorem~\ref
{thm-Bondareva-Shapley-megszamjatek-Schmeidler}.  The aforegiven notions
and concepts are utilized in the poof, and it will be seen that its main
ingredience is the use of compactness.

\begin {proof} [Proof of Theorem~\ref{thm-Bondareva-Shapley-megszamjatek-Schmeidler}]
Assume that the given coalition function $v\colon \mathcal A \to \mathbb
R$ is bounded below by~$L$.  We are to show that $\bacore(v) \neq
\emptyset$ if and only if the given game is balanced.  The
``only if'' part is obvious.  Assume that a $\mu \in \bacore(v)$ and let
$\mathcal S = \{S_1, S_2, \ldots, \allowbreak S_r\} \subseteq \mathcal
A$ be a balanced collection of coalitions, so that \eqref
{eq-Schmeidler-balanced-set-of-coalitions}~holds for some non-negative
balancing weights $\lambda_1$,~$\lambda_2$,~\dots,~$\lambda_r$.  Then
$\sum_{p=1}^r \lambda_p v(S_p) \leq \sum_{p=1}^r \lambda_p \mu(S_p) =
\mu(N) = v(N)$, so \eqref {eq-Schmeidler-balanced-game}~is satisfied,
and the game is balanced.  It remains to prove the ``if''
part.

Pick up any sets $S_0, S_1, \ldots, \allowbreak S_n \in \mathcal A$.
Our purpose is to show that $\bigcap_{j=0}^n F_{S_j} \neq \emptyset$.
We can assume w.l.o.g.\ that the sets $S_0$,~\dots,~$S_n$ are distinct
with $S_0 = \emptyset$ and $S_n = N$, and that the collection $\{S_0,
\ldots, \allowbreak S_n\} \subseteq \mathcal A$ is a field of sets.
(Roughly speaking, the more sets we pick up, the smaller the
intersection $\bigcap_{j=0}^n F_{S_j}$ is.  Having to show the
intersection is non-empty anyway, we can include the empty and the grand
coalition among the sets.  Moreover, we can add further sets
from~$\mathcal A$ so that the collection $\{S_0, \ldots, S_n\}$ becomes
a finite field of sets.  

We can also assume w.l.o.g.\ that
$S_1$,~\dots,~$S_{n'}$ are all the atoms of the field; that is, they are
all the minimal elements in the collection $\{S_1, \ldots, S_n\}$.
Obviously, the atoms $S_1$,~\dots,~$S_{n'}$ are pairwise disjoint, and
it holds $n = 2^{\smash{n'}\vphantom{n}} -\nobreak 1$.

Now, the sets $\emptyset = S_0$, $S_1$,~\dots,~$S_n$ being fixed, we
apply balancedness to the sets $S_1$,~\dots,~$S_n$:
\begin {equation}
\label {megszamjatek-Schmeidler-balancedness}
\vcenter {
	\normalbaselines
	\openup\jot
	\mathsurround=0pt
	\setbox0 = \hbox {%
		\eqref {megszamjatek-Schmeidler-balancedness}%
		$\mkern36mu$%
		\eqref {megszamjatek-Schmeidler-balancedness}}
	\dimen0 = \displaywidth
	\advance \dimen0 by -\wd0
	\ialign {\hbox to \dimen0 {$\strut #$}\cr
\forall \lambda_1, \ldots, \lambda_n \geq 0\colon\,\,\,
\lambda_1 \chi_{S_1} + \cdots + \lambda_n \chi_{S_n} = \chi_N
\hfil \cr \hfil
\mkern\thinmuskip \mkern\thickmuskip \Longrightarrow \mkern\thickmuskip
\lambda_1 v(S_1) + \cdots + \lambda_n v(S_n) \leq v(N) \,\hbox{.} \cr } }
\end {equation}
It shall follow hence that the system of relations
\begin {equation}
\label {megszamjatek-Schmeidler-core}
\eqarraywithdef {
	$\mu(#\hfil)$& \hfil ${}#{}$\hfil & $v(#\hfil)$& $\,\hbox{#}$\hfil \cr
\hfil N & =    & \hfil N & ,\cr
S_1     & \geq & S_1     & ,\cr
\multispan4\dotfill\cr
S_n     & \geq & S_n     &  }
\end {equation}
has a solution $\mu \in \ba(\mathcal A)$ such that
$\mathopen\|\mu\mathclose\| \leq R$.  To see that, we apply the
Bondareva-Shapley Theorem for finite games (Theorem~\ref
{thm-Bondareva-Shapley}).

Consider a new finite game $v'\colon \mathcal P(N') \to \mathbb R$ with
the set of the players $N' = \{1, \ldots, n'\}$.  Define the game as
follows.  Recall first that the collection $\{S_0, \ldots, S_n\}$ is a
field of sets and that $S_1$,~\dots,~$S_{n'}$ are all its atoms, which
are pairwise disjoint.  Now, for an $S' \subseteq N'$, let $S =
\bigcup_{i'\in S'} S_{i'}$, notice $S \in \{S_0, \ldots, S_n\}$, and put
$v'(S') = v(S)$.  The new finite game $v'$~has been defined thus.

Now, condition \eqref
{megszamjatek-Schmeidler-balancedness}~equivalently says that the new
game $v'$~is balanced.  By the Bondareva-Shapley Theorem (Theorem~\ref
{thm-Bondareva-Shapley}), its core is non-empty: there exist $a_1,
\ldots, \allowbreak a_{n'} \in \mathbb R$ such that
$\sum_{i'=1}^{\smash{n'}\vphantom{n}} a_{i'} = v'(N')$ and $\sum_{i'\in
S'} a_{i'} \geq v'(S')$ for any $S' \subseteq N'$.

The atoms $S_1$,~\dots,~$S_{n'}$ being non-empty sets, there exist
elements $x_{i'} \in S_{i'}$ for $i' = 1$,~\dots,~$n'$.  Consider the
measure
$$
\mu = a_1 \delta_{x_1} + \cdots + a_{n'} \delta_{x_{n'}} \,\hbox{,}
$$
where $\delta_{x_{i'}}$~is the Dirac measure concentrated at~$x_{i'}$.
We have $\mu(N) = \mu(S_1 \cup \cdots \cup S_{n'}) = a_1 + \cdots +
a_{n'} = v(N)$.  For any $j = 1$,~\dots,~$n$, let $S'_j = \{\, i' \in N'
: S_{i'} \subseteq S_j \,\}$.  Then $S_j = \bigcup_{i'\in\smash{S'_j}}
S_{i'}$, and $\mu(S_j) = \sum_{i'\in\smash{S'_j}} a_{i'} \geq v'(S'_j) =
v(S_j)$.  We have shown thus that the $\mu$~is a solution to the
system of inequalities~\eqref{megszamjatek-Schmeidler-core}.

Finally, let us calculate the norm~$\mathopen\|\mu\mathclose\|$ of the
solution, see~\eqref {eq-mu-norm}.  For a $T \in \mathcal A$, we observe
that
$$
\mu(T) = \sum_{\substack{ \smash{i'}=1 \\ x_{i'}\in T}}^{\smash{n'}\vphantom{n}} a_{i'} \,\hbox{.}
$$
Given pairwise disjoint sets $T_1, \ldots, T_s \in \mathcal A$ such that
$N = T_1 \cup \cdots \cup T_s$, and recalling
$\sum_{i'=1}^{\smash{n'}\vphantom{n}} a_{i'} = v(N)$, we have
$$
\eqarraywithdef {\hfil $\displaystyle #$& $\displaystyle {}#$\hfil \cr
\sum_{q=1}^s \bigl|\mu(T_q)\bigr| &=
\sum_{q=1}^s \,\,\, \biggl| \,
\sum_{\substack{ \smash{i'}=1 \\ \hidewidth x_{i'}\in T_q \hidewidth}}^{\smash{n'}\vphantom{n}}
a_{i'} \biggr| \leq
\sum_{q=1}^s
\sum_{\substack{ \smash{i'}=1 \\ \hidewidth x_{i'}\in T_q \hidewidth}}^{\smash{n'}\vphantom{n}}
\mathopen|a_{i'}\mathclose| \cr &=
\sum_{\smash{i'}=1}^{\smash{n'}\vphantom{n}} \mathopen|a_{i'}\mathclose| =
\sum_{\substack{ \smash{i'}=1 \\ \hidewidth a_{i'}\geq0 \hidewidth}}^{\smash{n'}\vphantom{n}} a_{i'} -
\sum_{\substack{ \smash{i'}=1 \\ \hidewidth a_{i'}  < 0 \hidewidth}}^{\smash{n'}\vphantom{n}} a_{i'} \cr &=
v(N) - 2
\sum_{\substack{ \smash{i'}=1 \\ \hidewidth a_{i'}  < 0 \hidewidth}}^{\smash{n'}\vphantom{n}} a_{i'} \leq
v(N) - 2 v\biggl(
\bigcup_{\substack{ \smash{i'}=1 \\ \hidewidth a_{i'}<0 \hidewidth}}^{\smash{n'}\vphantom{n}} S_{i'}
\biggr) \cr &\leq
v(N) - 2L = R \,\hbox{.} }
$$
It follows that $\mathopen\|\mu\mathclose\| \leq R$.  To conclude, we
have a $\mu \in \ba(\mathcal A)$ such that it is a solution
to~\eqref{megszamjatek-Schmeidler-core} and $\mathopen\|\mu\mathclose\|
\leq R$, which means $\mu \in \bigcap_{j=1}^n F_{S_j}$.  Since $F_{S_j}
\subseteq F_\emptyset$ for $j = 1$,~\dots,~$n$, it holds $\mu \in
\bigcap_{j=0}^n F_{S_j}$.  We have shown thus that the system
$\{F_S\}_{S\in\mathcal A}$ is centered.  As the closed $R$-ball $B_R$~is
weakly* compact, we have $\bacore(v) = \bigcap_{S\in\mathcal A} F_S \neq
\emptyset$.
\end{proof}


The following example demonstrates that Theorem~\ref
{thm-Bondareva-Shapley-megszamjatek-Schmeidler} cannot be generalized
further.  It presents an unbounded below game that is balanced, but its core is empty.

\begin{example}
\label {ex1}
Let the player set be $N = \mathbb N$, and let $\mathcal A = \{\, S
\subseteq N : S$~is finite or $N \setminus S$ is finite$\,\}$.  Consider
the game represented by the coalition function $v\colon \mathcal A
\to \mathbb R$ defined as follows: for any $S \in \mathcal A$, let
$$
v(S) =
	\begin{cases}
\hphantom{-} 1               & \text{if $S = \{1\}$,}\\
\hphantom{-} 1 + \frac{1}{n} & \text{if $S = \{1, n\}$ for $n = 2$,~$3$,~\dots,}\\
-\sum_{n\in T} \frac{1}{n}   & \text{if $S = N \setminus T$ for a finite $T \in \mathcal A$,}\\
\hphantom{-} 0               & \text{otherwise.}\\
	\end{cases}
$$
It is easy to see that this game is balanced.
Assuming that a $\mu \in \bacore(v)$, then $\mu\bigl(\{1\}\bigr) \geq
v\bigl(\{1\}\bigr) = 1$ and $\mu\bigl(N \setminus \{1\}\bigr) \geq
v\bigl(N \setminus \{1\}\bigr) = -1$.  Since $\mu\bigl(\{1\}\bigr) +
\mu\bigl(N \setminus \{1\}\bigr) = \mu(N) = v(N) = 0$, we have
$\mu\bigl(\{1\}\bigr) = 1$.  As $\mu\bigl(\{1,n\}\bigr) \geq
v\bigl(\{1,n\}\bigr) = 1 + 1/n$, it follows $\mu\bigl(\{n\}\bigr) \geq
1/n$ for all $n = 2$,~$3$,~\dots\spacefactor=\sfcode`. \space Summing
up, we have $\mu\bigl(\{1, \ldots, \allowbreak n\}\bigr) \geq \ln(n+1)$,
so $\mu \notin \ba(\mathcal{A})$ because $\mu$~is not bounded.  It
follows $\bacore(v) = \emptyset$.
\end{example}


The following ``limiting'' property of the ba-core is interesting.  It
is obtained as a corollary of Theorem~\ref
{thm-Bondareva-Shapley-megszamjatek-Schmeidler} by considering
the balancedness condition \eqref
{eq-Schmeidler-balanced-set-of-coalitions}~and~\eqref
{eq-Schmeidler-balanced-game}.


\begin {corollary}
\label {cor2}
Given a finite or infinite set\/~$N$ of the players and a field of
sets\/ $\mathcal A \subseteq \mathcal P(N)$ over\/~$N$, let\/ the
game represented by a coalition function\/ $v\colon \mathcal A \to
\mathbb R$ be bounded below.  For any\/ $\varepsilon > 0$, define the
coalition function\/ $v_\varepsilon\colon \mathcal A \to \mathbb R$ as
follows: let\/ $v_\varepsilon(N) = v(N) + \varepsilon$ and\/
$v_\varepsilon(S) = v(S)$ for all\/ $S \in \mathcal{A} \setminus \{N\}$.
Then: if\/ $\bacore(v_\varepsilon) \neq \emptyset$ for all\/
$\varepsilon > 0$, then\/ $\bacore(v) \neq \emptyset$.
\end {corollary}


Under the assumptions of Corollary~\ref{cor2}, the converse statement is
clear: if $\bacore(v) \neq \emptyset$, then $\bacore(v_\varepsilon) \neq
\emptyset$ for all $\varepsilon > 0$.  We thus conclude that a game
represented by the coalition function~$v$ is balanced if and
only if $\bacore(v_\varepsilon) \neq \emptyset$ for all $\varepsilon >
0$.

\section{Games with restricted cooperation}\label{sec6}

We now pay attention to games with restricted cooperation.  In
general, the cooperation is restricted whenever the collection $\mathcal
A \subseteq \mathcal P(N)$ of coalitions that can potentially emerge is
a proper subset of $\mathcal P(N)$.  In this sense, Theorem~\ref
{thm-Bondareva-Shapley-megszamjatek-Schmeidler} covers the case of
restricted cooperation, under the additional assumption that $\mathcal A
\subseteq \mathcal P(N)$ is a field of sets over~$N$, too.  Now, let
$\mathcal A' \subseteq \mathcal P(N)$ be the collection of all
coalitions that can potentially emerge; the collection $\mathcal
A'$~need not be a field of sets now.  Assume only $\emptyset, N \in
\mathcal A'$.  Then any coalition function $v'\colon \mathcal A' \to
\mathbb R$, such that $v(\emptyset) = 0$, represents a game with
restricted cooperation.

To introduce the concept of core of this game with restricted
cooperation, let $\mathcal A = \field(\mathcal A')$ be the field hull
of~$\mathcal A'$; that is, the minimal collection $\mathcal A \supseteq
\mathcal A'$ that is a field of sets over~$N$.  Then the core of a
game $v'$ with restricted cooperation is the set
$$
\eqarraywithdef {\hfil $# \mu($& $\hfil #\hfil$& $)\hfil #\hfil
	v'($& $\hfil #\hfil$& $)#\hfil $\cr
\bacore(v') = \bigl\{\, \mu \in \ba\bigl(\field(\mathcal A')\bigr) :
 & N & = & N & \,, \cr
 & S & \geq & S & \qquad \hbox {for all} \quad S \in \mathcal A' \setminus \{N\}
		\,\bigr\} \,\hbox{.} }
$$
We again ask whether $\bacore(v') \neq \emptyset$.

The following example presents a non-negative game with restricted cooperation which is balanced as by \eqref
{eq-Schmeidler-balanced-set-of-coalitions}~and~\eqref{eq-Schmeidler-balanced-game}, but only the feasible coalitions are considered, but its core is empty.

\begin {example}\label {exprc}
Let the player set be $N = \mathbb N$ and let $\mathcal A' =
\bigl\{\emptyset\bigr\} \cup \bigl\{\, \{1,i\} : i = 1,2, \allowbreak
3,\ldots\,\bigr\} \cup \bigl\{ N \setminus \{1\} \bigr\} \cup
\bigl\{N\bigr\}$.  Consider the game represented by the coalition
function $v'\colon \mathcal A' \to \mathbb R$ defined as follows: for
any $S \in \mathcal A$, let
$$
v'(S) =
	\begin{cases}
\hphantom{-} 2               & \text{if $S = \{1\}$,}\\
\hphantom{-} 2 + \frac{1}{n} & \text{if $S = \{1, n\}$ for $n = 2$,~$3$,~\dots,}\\
\hphantom{-} 0               & \text{if $S = N \setminus \{1\}$ or $S = \emptyset$,}\\
\hphantom{-} 1               & \text{if $S = N$.}\\
	\end{cases}
$$
Notice that this game is analogous to that presented in Example~\ref
{ex1}.  The field hull of~$\mathcal{A}'$ is $\mathcal{A} = \{\, S
\subseteq N : S$~is finite or $N \setminus S$ is finite$\,\}$.  The fact
that this game is balanced as by \eqref
{eq-Schmeidler-balanced-set-of-coalitions}~and~\eqref
{eq-Schmeidler-balanced-game}, but replace $\mathcal A$~and~$v$ with
$\mathcal A'$~and~$v'$, respectively, is clear.  To show that
$\bacore(v') = \emptyset$, it is enough to follow the arguments
presented in Example~\ref{ex1}.
\end{example}

Due to Example \ref{exprc}, we have to introduce a new notion of balancedness in the case of restricted cooperation. A game represented by a coalition function $v'$ defined on the class of feasible coalitions $\mathcal{A}'$ is bounded-balanced if there exists a bounded below balanced game $v$ defined  on $\mathcal{A}$, where $\mathcal{A}$ is the field hull of $\mathcal{A}'$ such that for every $S \in \mathcal{A}'$ it holds that $v (S) = v' (S)$. It is clear that if $\mathcal{A}'$ is a field and the game is bounded below (as in Theorem \ref{thm-Bondareva-Shapley-megszamjatek-Schmeidler}), then we get back to the notion of balancedness applied in Theorem \ref{thm-Bondareva-Shapley-megszamjatek-Schmeidler}. Moreover, notice that for finite games bounded-balancedness and balancedness by \citet{Faigle1989} are equivalent. 

The game in Example \ref{exprc} above has empty core because even if it is non-negative, none of its bounded below ``extensions" onto $\mathcal{A}$ is balanced and non of its balanced ``extensions" onto $\mathcal{A}$ is bounded below. 

%
%
%


Then the following theorem extends Theorem~\ref{thm-Bondareva-Shapley-megszamjatek-Schmeidler} to the class of games with restricted cooperation, hence it extends \citet{Faigle1989}.

\begin{theorem}\label {megszamjatek-Schmeidler-proposition}
Consider a coalition function\/ $v'\colon \mathcal A' \to \mathbb R$,
where\/ $\emptyset,N \in \mathcal A' \subseteq \mathcal P(N)$ and\/
$N$~is a finite or infinite set of the players. If\/ $v'$~is bounded below, then the\/ $\bacore(v') \neq \emptyset$ if and only if\/ $v'$~is bounded-balanced.
\end {theorem}

\begin{proof}
If $\mathcal{A}'$ is a field then we are back at Theorem~\ref{thm-Bondareva-Shapley-megszamjatek-Schmeidler}, hence nothing to do.

\bigskip

Suppose that  $\mathcal{A}'$ is not a field. The game $v'$ is bounded below, take any game  $v$ which makes $v'$ be bounded below. Then by Theorem \ref{thm-Bondareva-Shapley-megszamjatek-Schmeidler} $\bacore (v) \neq \emptyset$ if and only  $v$ is balanced. Since if $v$ is balanced then $v'$ is bounded-balanced by $v$, we get if $v'$~be bounded below, then the $\bacore{v'} \neq \emptyset$ if and only if $v'$ is bounded-balanced.
\end{proof}

Notice that the game $v'$~of Example~\ref {exprc} is not bounded-balanced, hence $\bacore(v') = \emptyset$.

\section* {Acknowledgements}

David Bartl acknowledges the support of the Czech Science Foundation
under grant number GA\v CR 21-03085S.  A part of this research was done
while he was visiting the Corvinus Institute for Advanced Studies; the
support of the CIAS during the stay is gratefully acknowledged.

Mikl\'os Pint\'er acknowledges the support by the Hungarian Scientific Research Fund under projects K 133882 and K 119930.


\end{document}